\documentclass[12pt]{amsart}
\usepackage{amsthm}
\usepackage{amsmath}
\usepackage{amssymb}
\usepackage{graphicx}
\usepackage{enumerate}
\allowdisplaybreaks
\numberwithin{equation}{section}
\numberwithin{figure}{section}

\theoremstyle{plain}
\newtheorem{theorem}{\sffamily Theorem}[section]
\newtheorem{proposition}[theorem]{\sffamily Proposition}
\newtheorem{lemma}[theorem]{\sffamily Lemma}
\newtheorem{corollary}[theorem]{\sffamily Corollary}
\newtheorem{example}[theorem]{\sffamily Example}
\newtheorem{remark}[theorem]{\sffamily Remark}
\newtheorem{definition}[theorem]{\sffamily Definition}
\newtheorem{conjecture}{\sffamily Conjecture}

\def\BET{\begin{theorem}}
\def\ENT{\end{theorem}}
\def\BEP{\begin{proposition}}
\def\ENP{\end{proposition}}
\def\BEL{\begin{lemma}}
\def\ENL{\end{lemma}}
\def\BEC{\begin{corollary}}
\def\ENC{\end{corollary}}
\def\BEE{\begin{example}\rm}
\def\ENE{\end{example}}
\def\BER{\begin{remark} \rm}
\def\ENR{\end{remark}}
\def\BED{\begin{definition} \rm}
\def\END{\end{definition}}
\def\BECJ{\begin{conjecture}}
\def\ENCJ{\end{conjecture}}

\def\bea{\begin{eqnarray}}
\def\eea{\end{eqnarray}}

\def\beas{\begin{eqnarray*}}
\def\eeas{\end{eqnarray*}}

\def\beq{\begin{equation}}
\def\eeq{\end{equation}}

\def\beal{\begin{align*}}

\def\eeal{ \end{align*} }

\def\roweq{\nonumber \\ &=& }
\def\rowleq{\nonumber \\  & \leq & }

\def\bfT{{\bf T}}

\def\bbC{{\mathbb C}}
\def\bbD{{\mathbb D}}

\def\bbN{{\mathbb N}}

\def\bbR{{\mathbb R}}

\def\cD{{\mathcal D}}

\begin{document}

%%%%%%%%%%%%%%%%%%%%%%%%%%%%%%%%%%%%%%%%%%%%%%%%%%%%%%%%%%%%%%%%%%%%%%%%%%%
%%    Please write the full title (english) in the title delimeter       %%
%%%%%%%%%%%%%%%%%%%%%%%%%%%%%%%%%%%%%%%%%%%%%%%%%%%%%%%%%%%%%%%%%%%%%%%%%%%

\title[Toeplitz operators on Bergman spaces]{
On generalized Toeplitz and little Hankel operators on Bergman spaces}

%%%%%%%%%%%%%%%%%%%%%%%%%%%%%%%%%%%%%%%%%%%%%%%%%%%%%%%%%%%%%%%%%%%%%%%%%%%
%% Please write the Mathematics Subject Clasification and keywords in    %%
%% the relative delimeters                                               %%
%% Put the acknowledgement before the bibliography if needed (see later) %%
%%%%%%%%%%%%%%%%%%%%%%%%%%%%%%%%%%%%%%%%%%%%%%%%%%%%%%%%%%%%%%%%%%%%%%%%%%%
%\def\classification{\falta}
%\def\keywords{\falta}

%%%%%%%%%%%%%%%%%%%%%%%%%%%%%%%%%%%%%%%%%%%%%%%%%%%%%%%%%%%%%%%%%%%%%%%%%%%
%%  Please write the authors                                             %%
%%  The authors adress must be after the bibliography (see later)        %%
%%  Please include a unique global support in the relative delimeter     %%
%%%%%%%%%%%%%%%%%%%%%%%%%%%%%%%%%%%%%%%%%%%%%%%%%%%%%%%%%%%%%%%%%%%%%%%%%%%

\author{Jari Taskinen}
\address{Department of Mathematics, University of Helsinki, 00014
Helsinki, Finland}
\email{jari.taskinen@helsinki.fi}

\author{Jani Virtanen}
\address{Department of Mathematics, 
 University of Reading,
 Whiteknights, P.O. Box 220,
 Reading RG6 6AX, UK}
\email{j.a.virtanen@reading.ac.uk}

%\thanks{JT was supported by the V\"ais\"al\"a Foundation of the Finnish Academy of Science %and Letters.}

%%%%%%%%%%%%%%%%%%%%%%%%%%%%%%%%%%%%%%%%%%%%%%%%%%%%%%%%%%%%%%%%%%%%%%%%%%%
%% Please write a short english abstract in the relative delimeters      %%
%%%%%%%%%%%%%%%%%%%%%%%%%%%%%%%%%%%%%%%%%%%%%%%%%%%%%%%%%%%%%%%%%%%%%%%%%%%
\begin{abstract}
We find a concrete integral formula for the class of generalized Toeplitz operators 
$T_a$ in Bergman spaces $A^p$, $1<p<\infty$, studied in an earlier work by the authors.
The result is  extended to little Hankel operators. 
We give an example of an $L^2$-symbol $a$ such that $T_{|a|} $ fails to be bounded in 
$A^2$, although  $T_a : A^2 \to A^2$ is  seen to be bounded by using the generalized 
definition. We also confirm  that the generalized definition coincides with the classical 
one  whenever the latter makes sense. 
\end{abstract}

\maketitle

%%%%%%%%%%%%%%%%%%%%%%%%%%%%%%%%%%%%%%%%%%%%%%%%%%%%%%%%%%%%%%%%%%%%%%%%%%%
%%                      THE TEXT                                         %%
%%      Use the following newtheorems if needed                          %%
%%      Please no change the style of these newtheorem commands          %%
%%      You may consider others newtheorem commands if needed            %%

\section{Introduction.}
\label{sec1}
Consider the space $L^p := ( L^p(\bbD, dA), \Vert \cdot \Vert_p)$,
where $1< p < \infty$ and $dA$ is the normalized area measure
on the unit disc $\bbD$ of the complex plane, and the Bergman space
$A^p$, which is the closed subspace of $L^p$ consisting of analytic functions. 
The Bergman projection $P$ is the orthogonal projection of
$L^2$ onto $A^2$, and it  has the integral representation
$$
    Pf(z) = \int_\bbD \frac{f(\zeta) }{ ( 1 - z \bar \zeta)^2 } dA (\zeta).
$$
It is also known to be a bounded projection of $L^p$ onto $A^p$ for every 
$1<p<\infty$.  For an integrable function $a : \bbD \to \bbC$ and, say, bounded
analytic functions $f$,  the Toeplitz operator 
$T_a$ with symbol $a$ is defined by 
\bea
T_af = P(af) = \int_\bbD \frac{a(\zeta) f(\zeta) }{ ( 1 - z \bar \zeta)^2 } dA (\zeta).
\label{0.0}
\eea
Since $P$ is bounded, it follows easily that $T_a$ extends to a 
bounded operator $A^p \to A^p$ for $1 < p < \infty$,  whenever $a$ 
is a bounded measurable function. %%A considerably more difficult 
The question of  the boundedness of $T_a$ on $A^p$ with unbounded symbols is a 
long-standing  problem. Examples of unbounded symbols inducing bounded Toeplitz 
operators can be easily constructed, since the behaviour of the symbol inside any 
compact subset of $\bbD$ is
not important for the boundedness of the operator. Also it is not difficult to find unbounded
symbols $a$ for which the integral in \eqref{0.0} converges, say, for all $f \in A^2$
but the operator is not bounded; see Section \ref{sec4} for an interesting example. 
We refer to the papers   \cite{E}, \cite{G},  \cite{G1},  \cite{Lu},  \cite{Lu2},
\cite{LT},  \cite{PTV1},  \cite{S},  \cite{SZ},  \cite{Su}, \cite{V},
\cite{V1}, \cite{Z1}, \cite{Zo} for  classical and recent results on the boundedness and
compactness of Toeplitz operators on Bergman spaces.

In the paper \cite{tasvir} we have given a generalized definition of 
Toeplitz operators, which we denote here $\bfT_a$. The definition takes 
efficiently  into account the possible cancellation phenomena of a
symbol. This leads to very weak  sufficient conditions for the boundedness of
Toeplitz operators. More precisely, in the reference it
was shown that  $\bfT_a$ is bounded under an averaging condition for the symbol itself
rather than for its modulus (the result is repeated and also extended 
to little Hankel operators  in Theorem \ref{th0.2}, below). However, the 
presentation of the result in \cite{tasvir} has some shortcomings and accordingly
the purpose of this paper is to make some improvements, which will be 
described in detail at the end of this section. 

The results of \cite{tasvir} show that 
cancellation phenomena may be essential in order to have a bounded operator $T_a$.
Here, we give an example which emphasizes this: in Section \ref{sec4}  we study the radial symbol $a \in L^2$, 
where $a(z) = |z|^{-1} (1-|z|)^{-1/4} \sin((1-|z|)^{-1}) $ for $|z| \geq 1/2$,  
and prove  that the operator $T_a$ is bounded in $A^p$, 
although $T_{|a|}$ is obviously not. Thus,  the boundedness of $T_a$ cannot be proven 
by conventional methods that only take into account the modulus of the symbol.
We can actually construct such a symbol in any given space $L^q$ with $q < \infty$.

Given  $ a \in L^1$, the little Hankel operator $h_a$ with symbol $a$ is
defined as 
$$
h_a f(z) = \int_\bbD \frac{a(\zeta) f(\zeta) }{ ( 1 - \bar z  \zeta)^2 } dA (\zeta)
$$
for $f \in A^p$ such that this integral converges. In this paper we make the
observation that the generalized definition of a Toeplitz operator and the  results of \cite{tasvir} can be extended to the little Hankel case as well. 
The results for $h_a$ are presented  in parallel with Toeplitz operators. 

As for the notation used in this paper, all function spaces are defined on the open unit disc
$\bbD$. In particular $H^\infty$ denotes the Hardy space of bounded analytic functions on 
$\bbD$.  If $0 < \rho < 1$, we denote $\bbD_\rho = \{ |z|\leq \rho \}$.
We also denote the standard weight by $W(z) = 1- |z|^2$, 
the kernel functions by $K_\lambda(z) =(1 - z \bar \lambda)^{-2}$
and $k_\lambda = K_\lambda / \Vert K_\lambda \Vert_2 = 
 W(\lambda) k_\lambda$,  and the M\"obius transform by
$\varphi_\lambda (z) = (\lambda -z  )/(1- z \bar \lambda)$, where 
$z, \lambda \in \bbD$. By $C$, $C'$ etc. we mean generic constants, the exact values of which 
may change from place to place.
We will deal with symbols $a$, which always at least belong to the space 
$L^1_{\rm loc}$ of locally integrable functions on $\bbD$. 
For other notation and definitions we refer to the book \cite{Z}.

Let us first describe briefly the  sufficient condition for the  boundedness
of generalized Toeplitz operators given in \cite{tasvir}.

\BED
\label{def0.3}
Denote by $\cD$ the family of the sets $D:= D(r, \theta )$ , where 
%%defined by
\bea
D=\{ \rho e^{i  \phi} \ | \
r \leq \rho \leq 1 -  \frac12  (1-r)\  ,
\ \theta \leq \phi \leq \theta  +    \pi  (1-r )  \}
\label{0.4}
\eea
for all $0 < r < 1$, $\theta \in [0,2 \pi]$.
We denote $ |D| := \int_{D} dA$ and, for  $\zeta = \rho e^{i\phi} \in D(r, \theta)$,
\bea
\hat a_D (\zeta) := \frac1{|D|} \int\limits_{r}^\rho \int\limits_\theta^\phi
a (\varrho e^{i\varphi}) \varrho d \varphi d \varrho  ,
\label{0.6}
\eea
where $a \in L^1_{\rm loc}$. In the following we will study symbols $a$ for which 
there exists a constant $C >0$ such that
\bea
|\hat a_D (\zeta) |\leq C  \label{0.7}
\eea
for all $D \in \cD$ and all $\zeta \in D$. 
\END

It turns out that one can proceed to a   generalized definition of bounded
Toeplitz operators just by using the condition \eqref{0.7}. However, for the 
proofs we  need to recall some more definitions from \cite{tasvir}. 
The countably many sets 
$D\big( 1 - 2^{- m + 1}, 2 \pi (\mu - 1 ) 2^{-m} \big) \in \cD$, 
where  $m \in \bbN , \mu = 1 , \ldots, 2^{-m}$, form a decomposition 
of the disc $\bbD$. We index these sets somehow into a family
$(D_n)_{n=1}^\infty$, so that every $D_n$ is of the form
\bea
D_n = &  \{ &  z = r e^{i\theta} \ | \
r_n < r \leq r_n' ,
\theta_n  < \theta \leq \theta_n'
\}
\label{2.2}
\eea
where, for some $m$ and $\mu$,
\bea
& & r_n  = 1 - 2^{- m + 1}
\ , r_n':= 1 -  2^{-m} , \ %\nonumber
\ \theta_n = \pi (\mu - 1 ) 2^{-m+1}\ ,  \  \theta_n' :=
 \pi \mu  2^{-m +1 }  . \label{2.2a}
\eea

Let $f \in A^p$. For all  $n = n(m,\mu)$ we write %%denote
\bea
& & F_n f (z) = \int\limits_{D_n}  \frac{a(\zeta) f(\zeta) }{
( 1 - z \bar \zeta)^2 } dA (\zeta)
\ \ , \ \  H_n f (z) = \int\limits_{D_n}  \frac{a(\zeta) f(\zeta) }{
( 1 - \bar z  \zeta)^2 } dA (\zeta) \ \ \ \ \forall \, z \in \bbD ,
\label{1.6}
\eea
so that $F_n$ can actually be considered as a conventional, bounded Toeplitz operator on
$A^p$; similarly for $H_n$.

Item $1^\circ$ of the following theorem is the main result Theorem 2.3
of \cite{tasvir}.
Also, $2^\circ$ is an immediate consequence of its proof: we leave to the 
reader the completely straightforward task to verify that the change 
$z \bar \zeta \to \bar z \zeta$ in the denominator does not affect the
proof.

\BET
\label{th0.2} Let $ 1 < p < \infty$ and 
assume that $a \in L^1_{\rm loc}$ satisfies the condition \eqref{0.7}. 
Then, the following hold true.

\noindent 
$1^\circ$. Given $f \in A^p$, the series $\sum_{n=1}^\infty F_nf (z)$ 
converges pointwise, absolutely for almost all $z \in \bbD$,
and the generalized Toeplitz operator $\bfT_a : A^p \to A^p$,
defined by
\bea
\bfT_a f (z) = \sum\limits_{n=1}^\infty  F_n f (z) \label{1.8}
\eea
is bounded
for all $1 < p < \infty$, and there is a constant $C$ such that
\beas
\Vert \bfT_a \Vert \leq C\sup\limits_{D \in \cD, \zeta \in D}
|\hat a_D (\zeta) | . %%\label{0.7a}
\eeas

$2^\circ$. For $f \in A^p$, the series $\sum_{n=1}^\infty H_nf(z)$ 
converges pointwise, absolutely,  for almost all $z \in \bbD$.
We define the generalized little Hankel operator by 
\bea
{\bf h}_a f (z) = \sum\limits_{n=1}^\infty  H_nf(z) \label{1.8h}
\eea
Then, ${\bf h}_a : A^p \to A^p$ is bounded
for all $1 < p < \infty$, and there is a constant $C$ such that
\beas
\Vert {\bf h}_a \Vert \leq C\sup\limits_{D \in \cD, \zeta \in D}
|\hat a_D (\zeta) | . %%\label{0.7h}
\eeas
\ENT

In this paper we improve Theorem \ref{th0.2} in the following ways.

$1^\circ$. The  definition \eqref{1.8} of a generalized  Toeplitz operator seems to depend on 
the geometry of a fixed decomposition \eqref{2.2} of the unit disc. (No doubt, other 
decompositions of $\bbD$, say with different choices of the points $r_n$ and  $\theta_n$, 
could be used as well, and it is not a priori clear, if the generalized operator  
defined in that way coincides with \eqref{1.8}. In fact,  an approach using Whitney decompositions with Euclidean rectangles for simply connected domains was presented in \cite{M}.) In this paper, formula 
\eqref{3.2}, we show  that the definition \eqref{1.8} coincides with a natural radial 
limit of conventional Toeplitz operators, and thus  the dependence of the definition on the 
decomposition  of the disc vanishes. 

$2^\circ$. It is not difficult to see that the generalized definition 
\eqref{1.8} of a Toeplitz operator coincides with the usual definition, whenever
the latter gives a bounded operator and condition \eqref{0.7} holds. This simple proof was 
omitted from \cite{tasvir}, but we present it here in Proposition \ref{th1.4}.

$3^\circ$. The terms $F_n$ in the series \eqref{1.8} are actually conventional, bounded
Toeplitz operators. In \cite{tasvir} it is only shown that the series \eqref{1.8} converges
in the very weak sense mentioned in  Theorem \ref{th0.2} above.  Here, we 
show in Theorem \ref{th1.3}
that the operator series $\sum_n F_n$  converges in the strong operator 
topology, and the same is true for the new limit representation \eqref{3.2}. 
Theorem \ref{th1.3} also contains an  immediate application of this result to transposed operators. 

$4^\circ$. The proof of Theorem 2.3 of \cite{tasvir} contains a small error: the 
inequality (3.8) of the citation is not true as such,
since the point $r_n' e^{i \theta_n'}$ there is actually on the boundary of the set $D_n$.
It is however not difficult to fix the flaw, and indeed in the course of the proof of 
Theorem \ref{th1.3} we do this by replacing the set $D_n$ by a bit larger set 
denoted  by $U_n$, see \eqref{3.3ar}.

\section{Main result.}
\label{sec2}

We now give a simplified expression  of the generalized Toeplitz operator 
$\bfT_a$, \eqref{1.8}, and also treat the little Hankel operator as well as the transposed
operators. 
Given  $a \in L_ {\rm loc}^1$  and $0 < \rho < 1$ we define 
the function $a_\rho : \bbD \to \bbC$ by $a_\rho (z) = a(z) $, if $|z| \leq \rho$ and 
$a_\rho(z) =0$ otherwise. It is plain that the  Toeplitz and little Hankel operators
\bea
T_{a_\rho} f (z) = \int\limits_{\bbD_\rho}  \frac{a(\zeta) f(\zeta) }{
( 1 - z \bar \zeta)^2 } dA (\zeta) \ \ , \ \
h_{a_\rho} f (z) = \int\limits_{\bbD_\rho}  \frac{a(\zeta) f(\zeta) }{
( 1 - \bar z  \zeta)^2 } dA (\zeta) \ \ , \ \ z \in \bbD , \label{3.1y}
\eea
are bounded $A^p \to A^p$.  

\BET 
\label{th1.3}
Let $ 1 < p < \infty$ and $1/p + 1/q =1$, and assume that $a \in L_ {\rm loc}^1$ and that 
\eqref{0.7} holds. Then, the generalized Toeplitz operators $\bfT_a : A^p \to A^p$ and little 
Hankel operators ${\bf h}_a : A^p \to L^p$,  defined in \eqref{1.8} and \eqref{1.8h}, respectively, can be written as 
\bea
\bfT_a f & = &  \lim\limits_{\rho \to 1} 
T_{a_\rho} f  \ , 
%\int\limits_{\bbD_\rho}  \frac{a(\zeta) f(\zeta) }{
%( 1 - z \bar \zeta)^2 } dA (\zeta)  \ , %%\mbox{and} \  
\label{3.2} \\
{\bf h}_a f&  = &  \lim\limits_{\rho \to 1} h_{a_\rho} f 
%\int\limits_{\bbD_\rho}  \frac{a(\zeta) f(\zeta) }{
%( 1 - \bar z  \zeta)^2 } dA (\zeta)
\label{3.2a}
\eea
for all $f \in A^p$. The limits converge with respect to  the strong operator topology
(SOT). %% and in particular pointwise for  almost all $z \in \bbD$.

%%Moreover, the operators 

Moreover, the transposed operators (with respect to the standard complex dual pairing)
$\bfT_a^*  : A^q \to A^q$ and  ${\bf h}_a^* : L^q  \to A^q$
can be written as 
\bea
\bfT_a^*  f  & = &  \lim\limits_{\rho \to 1} 
T_{\bar a_\rho} f \ , 
%\int\limits_{\bbD_\rho}  \frac{\overline{a} (\zeta) f(\zeta) }{
%( 1 - z \bar \zeta)^2 } dA (\zeta)  \ , %%\mbox{and} \  
\label{3.2b} \\
{\bf h}_a^*  g&  = &  \lim\limits_{\rho \to 1} 
h_{\bar a_\rho} f  
%\int\limits_{\bbD_\rho}  \frac{\overline{a}(\zeta) g(\zeta) }{
%( 1 - \bar z  \zeta)^2 } dA (\zeta) 
\label{3.2c}
\eea
for $ f \in A^q$ and $g \in L^q$, for  almost all $z \in \bbD$, and the limits 
here also converge in the SOT.
\ENT

{\bf Remark.}  In the course of the proof we also show that the sum in \eqref{1.8}
converges in the SOT  and thus improve the result of 
\cite{tasvir} also in this sense.  Of course, the limit
on the right hand side of \eqref{3.2} cannot in general 
converge in  the operator norm, since the operators $T_{a_\rho}$  are compact.

\bigskip

Proof. The proof will be given in a few steps. Moreover, 
we prove the statement \eqref{3.2} only for the Toeplitz operator, but
the reader is asked to observe the  necessary changes for the little Hankel case
\eqref{3.2a}. 

$(i)$ In the first step we review and strengthen the proof of 
Theorem 2.3 in \cite{tasvir} concerning the sum in \eqref{1.8}.
Let $f \in A^p$ be arbitrary.
%and at the same time we correct the small mistake in the inequality
%(3.8) of \cite{tasvir}. 

For all $n \in \bbN$ we define
the collection of all sets $D_\nu$ which touch the given $D_n$, 
more precisely,
\beas
\cD_n = \{ D_\nu \, : \, \nu \in \bbN, \ \overline D_\nu \cap \overline D_n \not= \emptyset \} .
\eeas
By the definition of the sets $D_n$, see \eqref{0.4}--\eqref{2.2a}, there exist
constants $N$, $M \in \bbN$ such that any set $\cD_n$ contains at most $N$ 
elements $D_\nu$ and  on the other hand, any set $D_\nu$ belongs to at most
$M$ sets $\cD_n$. Moreover, given $D_n$ and $w \in D_n$,  the subdomain 
\bea
{\textstyle \bigcup}_{D \in \cD_n} D  =: U_n
\label{3.3ar}
\eea
always contains a Euclidean disc $D(w,R)$ with center $w$ and radius 
$R=R(n)>0$ such that $|D(w,R)| \geq C |D_n|$ (use again the choice
of the sets $D_n$ to see this). 

We claim that for each $n$ and $w \in D_n$, 
\bea
|f( w) | \leq \frac{C}{|D_n|} \sum_{D \in \cD_n}
\int\limits_D |f(\zeta)| dA(\zeta) . \label{3.3a}
\eea
%% and  in particular for $w =  r_n'e^{u \theta_n'}$, see \eqref{2.2a}. 
To prove \eqref{3.3a}, let $D(w,R) \subset U_n$ be as above. Then, \eqref{3.3a} 
follows from the usual subharmonicity property for  $D(w,R)$:
\beas
|f( w) | \leq \frac{C}{|D(w,R)|}
\int\limits_{D(w,R)} |f(\zeta)| dA(\zeta)
\leq 
\frac{C'}{|D_n|} \int\limits_{U_n} |f(\zeta)| dA(\zeta) . %\label{3.3b}
\eeas
From now on we replace the incorrect inequality (3.8) of \cite{tasvir} by \eqref{3.3a}.

The proof of \cite{tasvir}, which uses the integration by parts -trick and the assumption 
\eqref{0.7}, yields the estimate %%CONSTANT MISSING ON THE RIGHT
\bea
& & |F_n f (z)| \leq C \sum_{D \in \cD_n} G_D (z) , \
\ \ \ \ \mbox{where}  \label{3.4}  \\
&&
 G_D (z)  = 
\int\limits_{D} \frac{|f(\zeta)| + |f'(\zeta)| \, W(\zeta) %%(1-|\zeta|^2)
 + |f''(\zeta)| \, W(\zeta)^2 %%(1-|\zeta|^2)^2
 }{|1 - z \bar \zeta|^2}
dA(\zeta) .
\nonumber
\eea
We observe by Theorem 4.28 of \cite{Z} that the function $g: = |f| + |f'| \, W 
%%(1-|\zeta|^2)
 + |f''| \, W^2$ in the integrand belongs to $L^p$. Following the argument in \cite{tasvir}, 
the positive term series 
\bea
\sum_{n=1}^\infty G_{D_n}(z)   \label{3.5}
\eea
converges for almost all $z$ and defines a function which belongs to $L^p$, since it 
it is pointwise bounded by the maximal Bergman projection $|P|$ of $g$.
%% the above mentioned $L^p$-function. 
Thus we see that also the series 
\bea
\sum_{n=1}^\infty \sum_{D \in \cD_n} G_{D} (z) \label{3.5a}
\eea
converges for almost all $z$, and the sum belongs to  $L^p$. This follows from the 
convergence of  \eqref{3.5}, since the terms of \eqref{3.5a} consist of the 
positive expressions $G_{D_n}$, and any single  $G_{D_n}$ can occur at most $MN$ times in
\eqref{3.5a}, by the definition of the numbers $N$ and
$M$.

%%As a consequence, 
By  \eqref{3.4}, the convergence of \eqref{3.5a} implies the absolute convergence of the  
series $\sum_n F_n f (z)$ a.e.. We claim that the operator sequence 
$(T^{(m)})_{m=1}^\infty $ defined by
\bea
T^{(m)} f = \sum_{n=1}^m F_n f
\eea
converges to $\bfT_a $ in the  SOT, as $m \to \infty$.
Indeed, given $f \in A^p$ and any $z \in \bbD$, the difference
\bea
\big| \bfT_a f (z) - T^{(m)} f (z) \big|  = \Big| \sum_{n > m} F_n f(z) \Big|  = 
\Big| \int\limits_{V_m}  \frac{f(\zeta)}{(1 - z \bar \zeta)^2} dA(\zeta) \Big| ,
\label{3.7b}
\eea
where $V_m =  \cup_{n >m} D_n$, has by \eqref{3.4} the upper bound
\bea
C \int\limits_{V_\mu}  \frac{g(\zeta)}{|1 - z \bar \zeta)^2|} dA(\zeta)
= C \int\limits_{\bbD}  
\frac{\chi_{V_\mu} (\zeta) g(\zeta)}{| 1 - z \bar \zeta |^2} dA(\zeta)
= C |P| ( \chi_{V_\mu}  g ) (z)  ; \label{3.7c}
\eea
here,  $\mu$ is some positive integer with $\mu \to \infty$ as $m \to \infty$,
and $\chi_{V_\mu}$ is the characteristic function of the set $V_\mu$.
But we have $\Vert  \chi_{V_\mu}  g \Vert_p \to 0$ as $\mu \to \infty$,
by Lebesgue's  dominated convergence theorem. 
Since $|P|$ is a bounded operator, there also holds 
$\Vert |P| ( \chi_{V_\mu}  g )  \Vert_p \to 0$ as $\mu \to \infty$. 
Combining this with the estimates \eqref{3.7b}--\eqref{3.7c} we get that 
$\Vert  \bfT_a f - T^{(m)} f \Vert_p \to 0 $ as $ m \to \infty$, which proves
the claim.

$(ii)$ We next consider the relation of the limit in \eqref{3.2} with the sum \eqref{1.8}. 

Let us fix  $n$ for a moment. Inspecting the proof of \cite{tasvir}
we see that given any 
$(\tilde r, \tilde \theta)$ such that $r_n <  \tilde r < r_n'$ and $\theta_n < \tilde \theta < \theta_n'$,
the expression  
\beas
G_n(z,\tilde r,\tilde \theta) :=   \int\limits_{r_n}^{\tilde r} 
\int\limits_{\theta_n}^{\tilde \theta}
\frac{a(\varrho e^{i \varphi} ) f(\varrho e^{i \varphi}) }{
( 1 - z \varrho e^{ - i \varphi})^2 } \varrho d \varrho d  \varphi 
%%\label{3.6}
\eeas 
has the same upper bound as $ F_n f  (z)$ in \eqref{3.4} (cf. \eqref{1.6}), namely
\bea
|G_n(z,\tilde r,\tilde \theta)|  \leq C \sum_{D \in \cD_n} G_D (z).  \label{3.8}
\eea
To see this one has to make the
straightforward changes to the upper limits of integrals in (3.6)--(3.11)
of \cite{tasvir} and also use \eqref{3.3a}. This is left to the reader as an easy task.

Given $\rho$, the integral in \eqref{3.2} can be written as
\bea
\int\limits_{\bbD_\rho}  \frac{a(\zeta) f(\zeta) }{
( 1 - z \bar \zeta)^2 } dA (\zeta) = \sum_{n=1}^m F_n f (z) + \sum_{n= m+1}^K G_n(z, \rho, \theta_n') \label{3.12}
\eea
for some integers $m$ and $K > m $, and moreover, $m \to \infty$ as
$\rho \to 1$. It is then obvious from the estimate \eqref{3.8} and
the convergence \eqref{3.5a} that  for almost all $z$, 
the limit in \eqref{3.2} must exist and, by \eqref{3.12},  it has to coincide
with $\sum_n F_n f(z) = \bfT_a f(z)$, \eqref{1.8}. 

Concerning the convergence in the  SOT, we use \eqref{3.8} and
\eqref{3.12} and the argument around \eqref{3.7b}--\eqref{3.7c} to estimate  
the difference
\bea
|\bfT_a f(z) - T_{a_\rho} f(z)| \leq C |P| (\chi_{V_\mu}g)(z)
\eea
where $\mu \to \infty$ as $\rho \to1$. Convergence in the  SOT follows in the same way as 
at the end of part $(i)$.

$(iii)$ Let us consider \eqref{3.2b}; let $f \in A^p$ and $g \in A^q$ be given. 
Denoting by $\langle \cdot , \cdot \rangle$ the standard complex dual paring of 
$A^p$ and $A^q$, we have 
\bea
& & \langle f, \bfT_a^* g \rangle = \langle \bfT_a f,  g \rangle
 = %\int\limits_{\bbD} \lim\limits_{\rho \to 1} 
%\int\limits_{\bbD_\rho}  \frac{a(\zeta) f(\zeta) }{
%( 1 - z \bar \zeta)^2 } dA (\zeta) \overline{g(z)} dA(z)
%\roweq 
\int\limits_\bbD  \bar g \lim\limits_{\rho \to 1} T_{a_\rho} f dA 
=  \lim\limits_{\rho \to 1} \int\limits_\bbD  \bar g \,  T_{a_\rho} f \, dA  , \label{3.80}
\eea
where the limit and the  integral could  be commuted because of the convergence of 
\eqref{3.2} in the SOT. 
Then, \eqref{3.80} equals
\beas
& &  \lim\limits_{\rho \to 1}  \int\limits_{\bbD}
\int\limits_{\bbD_\rho}  \frac{a(\zeta) f(\zeta) }{
( 1 - z \bar \zeta)^2 } dA (\zeta) \overline{g(z)} dA(z) 
\roweq 
\lim\limits_{\rho \to 1} 
\int\limits_{\bbD}  f(\zeta) a_\rho(\zeta)
\overline{ P g (\zeta)  }  dA(\zeta)
%% \frac{  g(z)} {( 1 - \bar z  \zeta)^2 } dA (z)}
%%\roweq
= \lim\limits_{\rho \to 1} \int\limits_{\bbD}  f(\zeta) a_\rho(\zeta)\overline{  g(\zeta )} dA (\zeta) 
\roweq
\lim\limits_{\rho \to 1} \int\limits_{\bbD} f(\zeta ) 
\overline{ P (\overline{ a_\rho} g)(\zeta)  } dA (\zeta) 
=   \lim\limits_{\rho \to 1} \int\limits_\bbD  f \, \overline{ T_{\bar a_\rho} g} \, dA
=   \int\limits_\bbD  f \,  \lim\limits_{\rho \to 1} \overline{ T_{\bar a_\rho} g} \, dA
 ,  %%\label{3.82}
\eeas
where at the end we used the fact that $\bar a$  obviously also satisfies condition \eqref{0.7} and the convergence of \eqref{3.2} in the SOT.

That the limit exist in the SOT follows from the treatment of the limit \eqref{3.2},
since $\bar a$  satisfies  \eqref{0.7}.
The proof of the little-Hankel case \eqref{3.2c} is similar, with obvious changes.
 \ \ $\Box$

\bigskip

\section{Concluding remarks.}
\label{sec4}

The following observation can be summarized  as saying  that $T_a f$ and $\bfT_a f $ 
coincide, whenever the former operator is bounded and condition \eqref{0.7} holds. 
%, since by the boundedness
%of $T_a$ we also have $T_a f_n \to T_a f $ in $A^p$. Also, $h_a$ and ${\bf h}_a$
%coincide in this sense.

\BEP
\label{th1.4} 
Let $1 < p < \infty$. 
Assume that $a \in L^1$, the integral \eqref{0.0} converges for all  $f \in A^p$
and  $T_a : A^p \to A^p$ is bounded; assume
moreover that \eqref{0.7} is satisfied  so that also $\bfT_a$ is bounded
in $A^p$.   Let 
$f \in A^p$ be arbitrary and then let $(f_n)_{n=1}^\infty \subset H^\infty$ 
be such that $f_n \to  f$ in $A^p$ as $n \to \infty$.  Then, 
$T_a f_n \to \bfT_a f $ in $A^p$, and, consequently, $T_a f = \bfT_af$ for all
$f \in A^p$.

The statement remains true for little Hankel operators, with $h_a$  replacing $T_a$  and
 ${\bf h}_a$ replacing $\bfT_a$.
 
%If in addition to the above assumptions, $a \in L^q$ for some $q $ with 
%$1/p \leq q < \infty $, then the defining integral \eqref{0.0} converges for 
%all $z$, the operator  $T_a : A^p \to A^p $ is bounded and we have 
%$T_a f = \bfT_af$ for all $f \in A^p$.
\ENP

Proof. Since $\bfT_a $ is a bounded operator $ A^p \to A^p$, we have 
$\bfT_a f_n \to \bfT_a f $ in $A^p$, and thus it is enough to show that 
$T_a g =  \bfT_a g $ for all $ g \in H^\infty$. But for such $g$,
the integral 
$$
   \int_\bbD \frac{a(\zeta) g(\zeta) }{ ( 1 - z \bar \zeta)^2 } dA (\zeta).
$$
converges, since $a \in L^1$ and the kernel function  $\zeta \mapsto
( 1 - z \bar \zeta)^2$ is bounded. Then  it is clear, see e.g. \cite{Ru}, Theorem 1.27, that 
$$
\sum_{n=1}^\infty \int_{D_n} \frac{a(\zeta) g(\zeta) }{ ( 1 - z \bar \zeta)^2 } dA (\zeta)   = \int_\bbD \frac{a(\zeta) g(\zeta) }{ ( 1 - z \bar \zeta)^2 } dA (\zeta)
= T_a g (z).
$$
This proves the result, since 
$$
\sum_{n=1}^N \int_{D_n} \frac{a(\zeta) g(\zeta) }{ ( 1 - z \bar \zeta)^2 } dA (\zeta)  
\to \bfT_a g(z) \ \ \ \mbox{as } \ N \to \infty 
$$
by what is mentioned around \eqref{1.8}.

The proof in the case of little Hankel operators is the same. 
 \ \ $\Box$

\bigskip

The sufficient condition \eqref{0.7} and the definitions \eqref{1.8}, \eqref{3.2} 
of Toeplitz operators are formulated for quite general locally integrable symbols, 
but the following example shows 
that the condition and the boundedness result are useful already in very simple, concrete
cases. 
A well known sufficient condition for the boundedness of $T_a$ is that
\bea
\sup\limits_{D \in \cD} M_a(D) < \infty \ \ \mbox{where} \
 M_a(D)  :=  %%\widehat{|a|} (D) = 
 \frac{1}{|D|} \int\limits_{D} |a|\, dA , \label{1.32}
\eea
and this condition is also necessary, if $\bbR \ni a(z) \geq 0$ for all $z \in \bbD$.
See \cite{Z}.

For every $0<b \leq 1/2$ we define the symbol
\begin{equation}
a_b(r e^{i \theta} ) := \left\{
\begin{array}{ll}
{\displaystyle \frac1{r (1-r)^b} \sin \frac1{1-r}\ , } & 
{\displaystyle \ \ r\geq \frac 12  } \\
{\displaystyle 1 \ , } & 
{\displaystyle \ \ r <  \frac 12 }
\end{array} \right. \label{5.30}
\end{equation}
which obviously belongs to $L^q$, if  $bq < 1$. Then, in particular, 
%%if e.g. $b = 1/4$,
 $a_{1/4} \in L^2$ and the defining integral formula of $T_{a_{1/4}}$ converges for
every $f \in A^2$. Obviously, the  defining formula of $T_{|a_{1/4}| }$ also converges for
every $f \in A^2$. However, we have the following result.

\BEP
\label{prop1.5}
$(i)$ The Toeplitz operator $T_{|a_b|}$ is not bounded in $A^p$ for any $1< p<\infty$
and  $0<b \leq 1/2$. 

\noindent $(ii)$ 
The Toeplitz operator $T_{a_b} $ is  bounded in $A^p$ for all $1< p<\infty$
and  $0<b \leq 1/2$. 
\ENP

Proof. Let us first deal with  $T_{|a_b|}$.  Given $0 < r < 1$ and any $\theta \in [0,2 
\pi]$, we consider the behaviour of $a_b$ in the set $D=D(r,\theta)$, see Definition
\ref{def0.3}. It is plain from the definition of $a_b$  and the elementary properties of 
the sinus  that for some universal constant $C >0$ we have 
\beas
|a_b( z )| \geq \frac14 (1-|z|)^{-b}  %%\label{5.32}
\eeas
in a subset of $D$ with area measure at least $C(1-r)^2$ (recall that $ |D|$ is proportional
to $(1-r)^2 $).
Then, of course $ M_a (D) \geq C' (1 -r)^{-b}$ for another constant $C' >0$, and
thus condition \eqref{1.32} cannot hold, and the operator $T_{|a_b|}$ is unbounded.

The symbol satisfies \eqref{0.7}, since 
given $D = D(1- 2 \delta , \theta)$ with a small enough $\delta$
and $\zeta = \rho e^{i \phi} \in D$,
we have, using the change of variable $y = 1/(1- \varrho)$ 
(so that $ d\varrho =  y^{-2} dy$) 
\beas
|D| | \hat a_D (\zeta) | & = &  \int\limits_\theta^\phi d\varphi 
\Big| \int\limits_{1 -2 \delta}^\rho 
\frac1{(1-\varrho)^b} \sin \frac1{1-\varrho}   d \varrho \Big| 
 = \pi \delta \Big| \int\limits_{1/(2 \delta)}^{1/(1-\rho)} 
y^{b-2} %\frac1{y^{1 + b } } 
\sin y \, dy \Big|  %\label{1.34}
\eeas
Let us divide the integration interval to subintervals
$J_n := [2 \pi n , 2 \pi(n+1 )]$, $n \in \bbN$. On $J_n$ we 
integrate as follows:   
\beas
& & \Big| \int\limits_{J_n } 
y^{b-2} %%\frac1{y^{1+ b}} 
\sin y \, dy \Big|
= \Big| \int\limits_{J_n} 
y^{b-2} %%\frac1{y^{1+ b}} 
\sin y \, dy  - (2\pi(n+1))^{b-2}%%\frac1{(2\pi(n+1))^{1+b}}
\int\limits_{J_n}  \sin y \, dy\Big|
\rowleq
 \int\limits_{J_n} \Big|
y^{b-2} %%\frac1{y^{1+b}}  
- (2\pi(n+1))^{b-2} %% \frac1{(2\pi(n+1))^{1+b}} 
\Big| \, dy \leq C n^{b - 3} .
 % \label{1.36}
\eeas
Hence, 
\bea
|D| | \hat a_D (\zeta) | & \leq &
\pi \delta \Big| \int\limits_{1/(2 \delta)}^{1/(1-\rho)} 
y^{b-2} \sin y \, dy \Big|
\leq 
C \delta \sum_{n= [1/(4 \pi \delta)]}^\infty n^{ b- 3 } \leq C' \delta^{3 - b} , \nonumber 
\eea
where $[x] $ denotes the integer part of a number $x \in \bbR$. Since $|D|$ is proportional 
to $\delta^2$,  the condition \eqref{0.7} holds true, and $T_{a_b}$ is bounded, by Theorem \ref{th0.2} and Proposition \ref{th1.4}.
\ \ $\Box$

\bigskip

{\bf Acknowledgements.} 
The authors wish to thank Grigori Rozenblum (G\"oteborg) for some personal communication
which initiated the investigation leading to  this work. The authors are also 
grateful for the  anonymous referees for remarks that helped to improve the 
results of this paper. 

The research of  JT was partially supported by the V\"ais\"al\"a Foundation of the Finnish Academy of Science and Letters. The research of JV and the visit of JT to the University of Reading were  supported by EPSRC grant EP/M024784/1.

\end{document}